\documentclass{amsart}
\usepackage[utf8]{inputenc}

\usepackage{amsmath, mathrsfs, amsfonts, mathtools, amsthm}
\usepackage{stmaryrd}
\usepackage{mathtools}
\usepackage{tikz-cd}
\usepackage{graphicx}
\usepackage{hyperref}
\usepackage[colorinlistoftodos]{todonotes}
\usepackage{faktor}
\newcommand{\dr}{\text{dr}}

\newcommand{\nd}{\text{dim}_{\text{nuc}}}

\newcommand{\K}{\mathcal{K}}
\newcommand{\T}{\mathcal{T}}
\newcommand{\Hil}{\mathcal{H}}

\author{Ruaridh Gardner and Aaron Tikuisis}	
\title[The nuclear dimension of extensions of $C(X)$ by $\mathcal K$]{The nuclear dimension of extensions of commutative C*-algebras by the compact operators}
\date{}

\newtheorem{lemma}{Lemma}
\newtheorem{theorem}[lemma]{Theorem}
\newtheorem{corollary}[lemma]{Corollary}

\newtheorem{definition}[lemma]{Definition}

\numberwithin{equation}{section}
\numberwithin{lemma}{section}

\renewcommand*{\thelemma}{\Alph{lemma}}

\thanks{Research partially supported by an NSERC discovery grant, held by AT.}

\begin{document}

\begin{abstract}
Generalizing the case of the Toeplitz algebra by Brake and Winter, we prove that the nuclear dimension of a C*-algebra extension of $C(X)$ by the compact operators is equal to the dimension of $X$.
\end{abstract}

\maketitle

Nuclear dimension is a numerical property of C*-algebras defined by Winter and Zacharias (\cite{WZ}) (inspired by an earlier notion by Kirchberg and Winter \cite{KW}).
Created by cleverly fusing the completely positive approximation property for C*-algebras together with Lebesgue's covering dimension, it is an intriguing concept that has played a prominent role in the classification of C*-algebras, particularly through its appearance in the Toms--Winter conjecture (see \cite[Section 5]{W18}).
It has received a lot of attention in recent years (see \cite{BBSTWW,CE,CETWW,CuntzToeplitz,DS,EGM,EM,EGLN,Evington,Geffen,HW,MS14,SWW15,W12,W16}, for example).

While many results establish finite nuclear dimension for classes of C*-algebras, determining the precise value is often a difficult problem.
The case of extensions is an important example of the difficulty.
Winter and Zacharias gave an upper bound for the nuclear dimension of an extension in terms of the values for the ideal and quotient, but it is easily seen to be non-optimal in the commutative case where classical techniques apply.
Answering a question posed in the paper defining nuclear dimension, Brake and Winter made a breakthrough in this direction, by computing the exact value (one) for the nuclear dimension of the Toeplitz algebra (\cite{BW}).

The Toeplitz algebra is, perhaps, the most well-known of the class of algebras featured in Brown, Douglas, and Fillmore's seminal extension theory (\cite{BDF}).
Many other interesting extensions of commutative C*-algebras by the algebra of compact operators were exposed and studied in the work of Brown, Douglas, and Fillmore and others; indeed, extensions of the form $0 \to \mathcal K \to E \to C(X) \to 0$ can be classified in terms of the $K$-homology of $X$ (\cite{KS}), where $X$ is a compact metrizable space.

In this article, we generalize Brake and Winter's result to these other extensions of a commutative C*-algebra $C(X)$ by the compact operators, showing that the nuclear dimension of any such extension agrees with the covering dimension of $X$:

\begin{theorem}\label{thm:Main}
Let $X$ be a compact metrizable space and let $E$ be an extension of $C(X)$ by $\K$:
\begin{equation} 0 \to \K \to E \to C(X) \to 0. \end{equation}
Then $\nd E = \dim(X)$.
\end{theorem}

Our approach makes use of many of the ideas from \cite{BW}.
Our main innovation arises in a step where Brake and Winter extend (up to a small perturbation) approximately order zero maps into matrix algebras, extending from the domain $C(S^1)$ to the domain $B(S^1$) (Borel functions).
This effectively allows a certain 1-dimensional piece of the overall completely positive approximation to become zero-dimensional.
In the present argument, replacing $S^1$ by an arbitrary metrizable compact space $X$, we need to be able to similarly extend approximately order zero maps to matrix algebras, from the domain $C(X)$ to the domain $B(X)$.

Brake and Winter use a correspondence between order zero maps from $C(S^1)$ and normal elements, then apply a classical result of Lin saying that approximately normal elements (corresponding to approximately order zero maps) can be perturbed to (exactly) normal elements -- to which Borel functional calculus applies.
In our more general setting, order zero maps from $C(X)$ correspond to $^*$-homomorphisms from $C_0((0,1]\times X)$, and we make use of a more refined result of Lin: a uniqueness theorem for approximately multiplicative maps from commutative C*-algebras to matrix algebras (\cite{Lin}).

However, Lin's uniqueness theorem has two crucial hypotheses: one $K$-theoretic and the other a faithfulness condition.
The $K$-theoretic condition holds trivially in our case, by the contractibility of $(0,1]\times X$.
The faithfulness condition is originally written in \cite{Lin} in terms of a uniform bound on measures of balls of a given radius, though we reformulate it as a faithful trace condition (see Theorem \ref{thm:LinReformulated}).
This condition means that the theorem cannot be applied automatically as is done in Brake and Winter's argument.
We require an additional step (Section \ref{sec:Quasicentral}), to ensure that we produce cpc order zero maps with an appropriate faithfulness condition in order to use Lin's result.
The flexibility of an absorbed trivial extension --- whose presence is obtained using classical Brown--Douglas--Fillmore theory -- allows us to realize this extra step.

In \cite{CuntzToeplitz}, the result of Brake and Winter is generalized in a different direction: by replacing the Toeplitz algebra with a Cuntz--Toeplitz algebra (or more generally, an extension of a Kirchberg algebra by the compact operators).
A further generalization can be found in \cite{Evington}, where the ideal of compact operators is replaced by a general stable AF algebra, and the quotient is a possibly non-simple $\mathcal O_\infty$-stable algebra.
Their key innovation is to apply a uniqueness theorem due to Gabe (\cite{Gabe}) for maps from $C_0((0,1],B)$, where $B$ is a separable nuclear $\mathcal O_\infty$-stable C*-algebra, similar to our use of Lin's uniqueness theorem.

The paper is organized as follows.
Section \ref{sec:Prelim} contains preliminaries about nuclear dimension and extensions.
We go into Lin's uniqueness theorem and formulate the relevant consequence of it in Section \ref{sec:Lin}.
Section \ref{sec:Quasicentral} is where we produce appropriate quasicentral approximate units to achieve appropriate cpc order zero map.
Finally, we prove Theorem \ref{thm:Main} in Section \ref{sec:Proof}.

The authors would like to thank Sam Evington, Abe Ng, and Stuart White for comments on an early version of the paper.

\renewcommand*{\thetheorem}{\roman{lemma}}
\setcounter{lemma}{0}
\numberwithin{lemma}{section}

\bigskip
\section{Preliminaries}
\label{sec:Prelim}

\begin{definition}[{\cite[Definition 3.1]{KW},\cite[Definition 2.1]{WZ}},\cite{WZCPC0}]
Let $A,B$ be C*-algebras.
A cpc map $\phi:A \to B$ is \emph{order zero} if it is orthogonality-preserving: $\phi(a)\phi(b)=0$ whenever $a,b \in A_+$ satisfy $ab=0$.
The \emph{nuclear dimension} of $A$ is the minimum number $d$ such that for any finite set $\mathcal F\subset A$ and any $\varepsilon>0$, there exist finite dimensional C*-algebras $F^{(0)},\dots, F^{(d)}$, together with cpc maps $\psi^{(i)}:A \to F^{(i)}$ and cpc order zero maps $\phi^{(i)}:F^{(i)}\to A$, such that
\begin{equation} \left\|\sum_{i=0}^d \phi^{(i)}(\psi^{(i)}(x))-x\right\| < \varepsilon, \end{equation}
for all $x \in \mathcal F$.

The \emph{decomposition rank} of $A$ is defined in the same way, but additionally asking that the sum $\phi^{(0)}+\cdots+\phi^{(d)}:F^{(0)}\oplus\cdots\oplus F^{(d)}\to A$ is cpc.

We write $\nd A$ and $\dr(A)$ for the nuclear dimension and decomposition rank (respectively) of $A$.
\end{definition}

We write $\mathcal K$ to denote the C*-algebra of compact operators on a separable infinite dimensional Hilbert space (such as $\ell^2(\mathbb N)$).

Let $X$ be a compact metrizable space.
An \emph{extension} of $C(X)$ by $\mathcal K$ is a C*-algebra $E$ fitting into a short exact sequence
\begin{equation} 0 \to \mathcal K \to E \to C(X) \to 0. \end{equation}
Extensions of $C(X)$ by $\mathcal K$ correspond to injective $^*$-homomorphisms from $C(X)$ to the Calkin algebra (one direction of this correspondence is described below; see also \cite[Definition 1.5 and Remark 1.6]{BDF77}, for example).
By the Choi--Effros lifting theorem,\footnote{In fact, a motivating special case of this theorem, due to Arveson, \cite{Arve}, suffices here.} one can lift such a $^*$-homomorphism to a ucp map $\phi:C(X)\to B(\mathcal H)$, and the property of being an injective $^*$-homomorphism modulo $\K$ translates into the conditions
\begin{enumerate}
\item $\phi(fg)-\phi(f)\phi(g) \in \K \text{ for all } f, g \in C(X),$
\item $\phi(f)\in \K \text{ if and only if } f=0$.
\end{enumerate}

\begin{definition}\label{ucp ext}
Given a ucp map $\phi:C(X)\to B(\mathcal H)$ that is an injective $^*$-homomorphism modulo $\K$ (i.e., satisfying (i) and (ii)) above, define
\begin{equation} \T_\phi:=\mathcal K+\phi(C(X)). \end{equation}
(This is the corresponding extension of $C(X)$ by $\mathcal K$).
\end{definition}

When $\phi:C(X)\to B(\mathcal H)$ is a genuine $^*$-homomorphism (still injective modulo $\K$), then $\mathcal T_\phi$ is called an \emph{trivial extension}.
We note that Theorem \ref{thm:Main} follows easily from known results in the case of a trivial extension, as follows: for such an extension, one can find a quasicentral approximate unit (see Section \ref{sec:Quasicentral}) for $(\K,\T_\phi)$ consisting of projections, and then by \cite[Proposition 6.2]{KW}, $\nd \T_\phi =\dr(\T_\phi)= \mathrm{dim}(X)$.

Given two ucp maps $\phi,\psi:C(X) \to B(\mathcal H)$ satisfying (i) and (ii) above, define $\phi\oplus\psi:C(X) \to M_2\otimes B(\mathcal H) \cong B(\mathcal H)$; this ucp map will satisfy the same conditions (i) and (ii), and therefore also defines an extension of $C(X)$ by $\K$.

\begin{theorem}[Brown--Douglas--Fillmore]
\label{thm:BDF} 
Let $\T_\phi,\T_\psi$ be extensions of $C(X)$ by $\K$, with $\T_\psi$ trivial.
Then $\T_\phi \cong \T_{\phi\oplus\psi}$.
\end{theorem}

\begin{proof}
This result follows easily from the work of \cite{BDF77}, as we explain now.
By \cite[Theorem 1.17]{BDF77}, the extensions corresponding to $\phi$ and $\phi\oplus\psi$ are equivalent (in the sense of \cite[Definition 1.1]{BDF77}).
By \cite[Remark 1.6(i)]{BDF77}, this implies that $\T_\phi \cong \T_{\phi\oplus\psi}$.
\end{proof}

\bigskip
\section{Correcting approximately multiplicative maps}
\label{sec:Lin}

We make use of the following ``approximate uniqueness'' result of Lin; we explain the notation following the statement.

\begin{theorem}[{\cite[Theorem 2.10]{Lin}}]\label{thm:Lin}
Let $Y$ be a compact metric space, let $\varepsilon>0$ and let $\mathcal{F} \subset C(Y)$ be a finite subset.
There exists $\eta>0$ satisfying the following:  for any $\sigma >0$, there exist $\gamma > 0$, $\delta >0$, a finite subset $\mathcal{G} \subset C(Y)$, a finite subset $\Hil \subset C(Y)_{\text{s.a.}}$, and a finite subset $\mathcal{P} \subset \underline{K}(C(Y))$ satisfying the following:

For any ($\mathcal G,\delta$)-multiplicative ucp maps $\varphi, \psi:C(Y) \to M_n$ (for some integer $n \ge 1$) for which
\begin{align} 
\label{eq:lin-partially-defined}
\varphi_*|_{\mathcal{P}} = \psi_*|_{\mathcal{P}}, \\
\notag
\mu_{\tau\varphi}(O_r) \ge \sigma
\end{align}
for all open balls $O_r$ of radius $r \ge \eta$, and 
\begin{equation}
\label{eq:lin-traces-approx-agree}
 \left| \tau\varphi(a)- \tau\psi(a)\right| < \gamma \text{ for all } a \in \Hil, \end{equation}
where $\tau$ is the tracial state of $M_n$, there is a unitary $u \text{ in } M_n$ such that 
\begin{equation} \|\varphi(f) - Ad_u\psi(f)\| < \varepsilon \text{ for all } f \in \mathcal{F}.\end{equation}
\end{theorem}

The object $\underline{K}(A)$ is the \textit{`total $K$-theory'} of $A$, the direct sum of the $K$-theory of $A$ and its $K$-theory with coefficients, $K_i(A,\mathbb Z_n):= K_i(A\otimes C_n)$, where $C_n$ is a commutative C*-algebra such that $K_0(C_n)=\mathbb Z_n$ and $K_1(C_n)=0$ (see \cite{Sch}).

We say that a map between C*-algebras $\varphi:A \to B$ is \emph{$(\mathcal G,\delta)$-multiplicative} (where $\mathcal G \subset A$ and $\delta>0$) if
\begin{equation}
\|\varphi(xy)-\varphi(x)\varphi(y)\|<\delta, \quad x,y \in \mathcal G.
\end{equation}

The notation in \eqref{eq:lin-partially-defined} is a sort of shorthand.
Implicitly, one asks for (finitely many) representatives of appropriate unitaries and projections in $C(Y) \otimes M_m$ and $(C(Y) \otimes C_n)^\sim \otimes M_m$ in order to capture all the classes in the finite set $\mathcal P$; call all these elements $\mathcal F_0$.
Moreover, it is implicit that by making $\delta$ small enough and $\mathcal G$ large enough, each projection (resp. unitary) in $\mathcal F_0$ is sent close enough to a projection (resp. unitary) so that it has a well-defined $K_0$-class (respectively $K_1$-class), by any ($\mathcal G,\delta$)-multiplicative map (and in particular, by both $\varphi$ and $\psi$).
Then \eqref{eq:lin-partially-defined} is saying that the $K_0$-class or $K_1$-class of $\psi(a)$ (depending on whether $a$ is a projection or unitary) agrees with that of $\varphi(a)$, for $a \in \mathcal F_0$.

We wish to reformulate Lin's uniqueness theorem in the language of sequence algebras.
In preparation, we have the following definition and lemma.

\begin{definition}
Let $B_n$ be a sequence of C*-algebras such that $T(B_n)\neq \emptyset$ for all $n$.
Define $B_\infty:=\prod_n B_n/\bigoplus_n B_n$.
A \emph{limit trace} on $B_\infty$ is a trace $\tau \in T(B_\infty)$ of the form
\begin{equation} \tau((b_n)_{n\in\mathbb N}) := \lim_{n\to\omega} \tau_n(b_n), \end{equation}
where $\tau_n \in T(B_n)$ for each $n$ and $\omega$ is a free ultrafilter.
Define $T_\infty(B_\infty)$ to be the set of all limit traces on $B_\infty$.
\end{definition}

\begin{lemma}
\label{lem:SequenceAlgFacts}
Let $B_n$ be a sequence of unital C*-algebras such that $T(B_n)\neq \emptyset$ for all $n$.
Define $B_\infty:=\prod_n B_n/\bigoplus_n B_n$.
Let $A$ be a separable C*-algebra.
Let $\varphi_n,\psi_n:A \to B_n$ be a sequence of ucp maps, and let $\varphi,\psi:A \to B_\infty$ be the maps induced by these sequences.
Then the following hold.
\begin{enumerate}
\item $\varphi$ is a $^*$-homomorphism if and only if there is 
an increasing sequence $(\mathcal G_n)_{n\in\mathbb N}$ of finite subsets of $A$ such that $\bigcup_n \mathcal G_n$ is dense in $A$ and a sequence $(\delta_n)_{n\in\mathbb N}$ of positive numbers converging to $0$, such that $\varphi_n$ is $(\mathcal G_n,\delta_n)$-multiplicative for all $n$.
\item Suppose that $\varphi,\psi$ are $^*$-homomorphisms satisfying $\varphi_*(x)=\psi_*(x)$ for all $x \in \underline{K}(A)$.
Then for any finite set $\mathcal P \subset \underline{K}(A)$, there exists $n_0$ such that for all $n\geq n_0$, $\left(\varphi_n\right)_*|_{\mathcal{P}} = \left(\psi_n\right)_*|_{\mathcal{P}}$.
\item
$\tau\circ\varphi=\tau\circ\phi$ for all $\tau \in T_\infty(B_\infty)$ if and only if for every finite set $\mathcal H\subset A$ and $\gamma>0$, there exists $n_0$ such that for all $n\geq n_0$, $\left| \tau(\varphi_n(a))- \tau(\psi_n(a))\right| < \gamma$ for all $a \in \mathcal H$ and $\tau \in T(B_n)$.
\item
$\varphi$ is unitarily equivalent to $\psi$ if and only if for every finite subset $\mathcal F\subset A$ and $\varepsilon>0$, there exists $n_0$ such that, for all $n\geq n_0$, there is a unitary $u \in A$ such that $\|\varphi_n(a) - Ad_u(\psi_n(a))\| < \varepsilon$ for all $a \in \mathcal{F}$.
\item Suppose that $B_n$ has a unique trace, $\tau_n$, for all $n$ and $A=C(Y)$ for some compact metric space $Y$.
Let $\mu_n$ be the probability measure on $Y$ associated to the trace $\tau_n\circ\varphi_n \in T(C(Y))$.
Then $\inf_{\tau \in T_\infty(B_\infty)} \tau\circ\varphi(a)>0$ for all $a \in A_+\setminus\{0\}$ if and only if, for every $\eta>0$ there exists $\sigma>0$ and $n_0$ such that for all $n\geq n_0$ and every open ball $O_r$ in $Y$ of radius $r\geq \eta$, $\mu_n(O_r)>\sigma$.
\end{enumerate}
\end{lemma}

\begin{proof}
(i):
$\Rightarrow$: Take any increasing sequence $(\mathcal G_k)_{k=1}^\infty$ of finite subsets of $A$ such that $\mathcal G_0=\emptyset$ and $\bigcup_k \mathcal G_k$ is dense in $A$, and take any sequence $(\delta_k)_{k\in\mathbb N}$ of positive numbers converging to $0$.
Then for each $k$, we must have that $\varphi_n$ is $(\mathcal G_k,\delta_k)$-multiplicative, for $n$ sufficiently large.
We may therefore find a sequence $(k_n)_{n\in\mathbb N}$ converging (possibly slowly) to $\infty$, such that each $\varphi_n$ is $(\mathcal G_{k_n},\delta_{k_n})$-multiplicative.

$\Leftarrow$: This is immediate.

(ii):
It suffices to show that, for any projection (respectively unitary) $a$ in $A \otimes M_m$ or $(A \otimes C_n)^\sim \otimes M_m$ (as appropriate), the $K_i$-class of $\varphi_n(a)$ and $\psi_n(a)$ agree, for $n$ sufficiently large.
By replacing $A$ with $A \otimes M_m$ or $(A \otimes C_n)^\sim\otimes M_m$ if necessary, we may assume $a \in A$.
Since the $K_i$-class of $\varphi(a)$ and $\psi(a)$ are equal, by possibly enlarging the matrix size, we may assume that they are homotopic.
In particular, there exist projections (resp.\ unitaries) $b_1,\dots,b_k\in B_\infty$, such that $b_1=\varphi(a)$, $b_k=\psi(a)$, and $\|b_i-b_{i+1}\|$ is small, for all $i=1,\dots,k-1$.
By lifting these, there is $n_0$ such that for $n\geq n_0$, there exist approximate projections (resp.\ unitaries) $\hat b_1,\dots, \hat b_k$ in $B_n$ such that $\|\hat b_i-\hat b_{i+1}\|$ is small, for $i=1,\dots,k-1$, with $\hat b_1=\varphi_n(a)$ and $\hat b_k=\psi_n(a)$.
Consequently, $\varphi_n(a)$ and $\psi_n(a)$ have the same $K_i$-class, for $n\geq n_0$.

(iii):
$\Rightarrow$: Suppose $\tau\circ\varphi=\tau\circ\psi$ and, for a contradiction, there exists a finite set $\mathcal H\subset A$ and $\gamma>0$ and infinitely many $n$ for which there exist $a \in \Hil$ and $\tau_n \in T(B_n)$ satisfying $\left| \tau_n\varphi_n(a)- \tau_n\psi_n(a)\right| \geq \gamma$.
By passing to a smaller but still infinite set of these $n$, we may assume that there is a single $a \in \mathcal H$ for which $\left| \tau_n\varphi_n(a)- \tau_n\psi_n(a)\right| \geq \gamma$ for infinitely many $n$.
Taking an ultrafilter concentrated on these $n$, and using it with these $\tau_n$ to get a limit trace $\tau \in T_\infty(B_\infty)$, we see that $|\tau(\varphi(a))-\tau(\psi(a))|\geq \gamma>0$, a contradiction.

$\Leftarrow$: Fix $a \in A$.
By the hypothesis, $\lim_{n\to\infty} |\tau_n(\varphi_n(a))-\tau_n(\psi_n(a))|=0$ for any choice of traces $\tau_n \in T(B_n)$.
From this it is clear that $\tau\circ\varphi(a)=\tau\circ\psi(a)$ for all $\tau \in T_\infty(B_\infty)$.

(iv): This is entirely standard.

(v):
$\Rightarrow$: Given $\eta>0$, by the Lebesgue's Number Theorem, there is an open cover $U_1,\dots,U_k$ of $Y$ such that every open ball $O_r$ of radius $r\geq \eta$ contains one of the sets $U_i$.
Let $f_i\in C(Y)$ be a nonzero positive contractive function supported on $U_i$.
Define 
\[ \sigma:=\frac12\inf_{\stackrel{i=1,\dots,k}{\tau \in T_\infty(B_\infty)}} \tau\circ\varphi(f_i); \]
 this is positive by hypothesis.
By the argument in (iii), there exists $n_0$ such that, for all $n\geq n_0$ and all $i=1,\dots,k$, $\tau_n(\varphi_n(f_i))>\sigma$.
Given an open ball $O_r$ of radius $r\geq \eta$, it contains an open set $U_i$ for some $i$, and so
\begin{equation} \mu_n(O_r)\geq \tau_n(\varphi_n(f_i))>\sigma. \end{equation}

$\Leftarrow$: Fix $a \in A_+\setminus \{0\}$.
Let $\eta,\delta>0$ be such that there is an open ball $O_\eta$ of radius $\eta$ such that $a(x)\geq \delta$, for all $x \in O_\eta$.
Let $\sigma$ and $n_0$ be as in the hypothesis; then we have
\begin{equation} \tau_n(\varphi_n(a))\geq \delta\mu_n(O_\eta)>\delta\sigma, \end{equation}
for all $n\geq n_0$, and therefore $\tau(\varphi(a))\geq \delta\sigma>0$ for all $\tau \in T_\infty(B_\infty)$.
\end{proof}

Here is the recasting of Lin's uniqueness theorem;\footnote{To derive Lin's theorem from this statement, one needs a converse to Lemma \ref{lem:SequenceAlgFacts} (ii), which could be established using $K$-theoretic regularity properties of matrix algebras; since we do not need this, we do not pursue it further here.}
by avoiding the use of a metric on the space $X$, it allows us to also generalize it to the non-unital case.

\begin{theorem}
\label{thm:LinReformulated}
Let $(m_n)_{n\in\mathbb N}$ be a sequence of natural numbers and set $B_\infty:=\prod_n M_{m_n}/\bigoplus_n M_{m_n}$.
Let $Y$ be a second countable locally compact Hausdorff space, and let $\varphi,\psi:C_0(Y) \to B_\infty$ be $^*$-homomorphisms, such that $\tau\circ\varphi,\tau\circ\psi$ are faithful traces on $C_0(Y)$ for every $\tau \in T_\infty(B_\infty)$.
Then $\varphi,\psi$ are unitarily equivalent if and only if $\underline{K}(\varphi)=\underline{K}(\psi)$ and $\tau\circ\varphi=\tau\circ\psi$ for every $\tau \in T_\infty(B_\infty)$.
\end{theorem}

\begin{proof}
The forward implication is immediate, so let us consider the reverse: we assume that $\varphi,\psi$ agree on $\underline{K}$ and on traces.
If $Y$ is not compact, then the hypotheses ensure that the unitizations $\varphi^\sim,\psi^\sim:C_0(Y)^\sim \to B_\infty$ also agree on $\underline{K}$ and traces, and satisfy the faithful trace hypothesis (since the extension of a faithful trace on $C_0(Y)$ to $C_0(Y)^\sim$ is faithful).
Thus, by possibly replacing $Y$ by its one-point compactification, we may assume that $Y$ is compact.

Let us fix $a \in C(Y)_+\setminus \{0\}$ and let us show that $\inf_{\tau \in T_\infty(B_\infty)} \tau(\varphi(a))>0$.
Write $\varphi(a)=(b_n)_{n\in\mathbb N}$.
Suppose for a contradiction that there are traces $\tau_k \in T_\infty(B_\infty)$ such that $\tau_k(\varphi(a))\to 0$.
Then, since each $\tau_k$ is a limit trace, there exists $n_k$ arbitrarily large such that $\tau_{M_{m_{n_k}}}(b_{n_k})<\tau_k(\varphi(a))+1/k$.
We may thus choose $n_1 < n_2 < \cdots$.
If we let $\omega$ be any ultrafilter concentrated on $\{n_1,n_2,\dots\}$ and define $\tau_\omega\in T_\infty(B_\infty)$ to be the corresponding limit trace, it follows that $\tau_\omega(\varphi(a))\leq \liminf_k \tau_{M_{m_{n_k}}}(b_{n_k}) = 0$, contradicting the assumption that $\tau_\omega \circ\varphi$ is a faithful trace.

By the Choi--Effros Lifting Theorem (\cite[Theorem C.3]{BO}), we may lift $\varphi,\psi$ to sequences of cpc maps $\varphi_n,\psi_n:A\to M_{m_n}$.
Let us explain how to use Theorem \ref{thm:Lin} to verify the condition in Lemma \ref{lem:SequenceAlgFacts}(iv), which implies that $\varphi,\psi$ are unitarily equivalent.
Let $\mathcal F$ be a finite subset of $C(Y)$ and let $\varepsilon>0$.
Let $\eta>0$ be as in Theorem \ref{thm:Lin} for this $\mathcal F$ and $\varepsilon$.
By Lemma \ref{lem:SequenceAlgFacts}(v), let $\sigma>0$ and $n_0$ be such that for all $n\geq n_0$ and every open ball $O_r$ in $Y$ of radius $r\geq \eta$, $\mu_n(O_r)>\sigma$.
Let $\gamma,\delta>0$, $\mathcal G\subset C(Y)$, $\Hil \subset C(Y)_{\text{s.a.}}$ and $\mathcal P \subset \underline{K}(C(Y))$ be as in Theorem \ref{thm:Lin} for this $\sigma$.
Then by Lemma \ref{lem:SequenceAlgFacts}(i),(ii),(iii), we may possibly increase $n_0$ so that, for all $n\geq n_0$, $\varphi_n,\psi_n$ are ($\mathcal G,\delta)$-multiplicative and satisfy \eqref{eq:lin-partially-defined} and \eqref{eq:lin-traces-approx-agree}.
Thus for every $n\geq n_0$, by Theorem \ref{thm:Lin}, there is a unitary $u \in M_{m_n}$ such that $\|\varphi_n(f) - Ad_u\psi_n(f)\| < \varepsilon$ for all $f \in \mathcal{F}$; consequently by Lemma \ref{lem:SequenceAlgFacts}(iv), $\varphi,\psi$ are unitarily equivalent.
\end{proof}

The following is an easy tracial existence result for $^*$-homomorphisms from a commutative C*-algebra to a matrix algebra.
In the sequel, it will be important to know that all arbitrarily large matrix sizes can be used in the codomain.

\begin{lemma}
\label{lem:TracesClose}
Let $Y$ be a locally compact Hausdorff space.
Then for any finite subset $\mathcal H \subset C_0(Y)$ and any $\gamma>0$ there exists $m_0$ such that, for any $\tau \in T(C_0(Y))$ and any $m\geq m_0$, 
there is a $^*$-homomorphism $\psi:C_0(Y) \to M_m$ such that
\begin{equation}
\label{eq:TracesClose}
 |\tau(f)-\tau_{M_m}(\psi(f))|<\gamma,\quad f \in \mathcal H. \end{equation}
\end{lemma}

\begin{proof}
Since $T(C_0(Y))$ is weak$^*$-compact, it suffices to prove that for each $\tau \in T(C_0(Y))$, there exists $m_0$ such that for any $m \geq m_0$ there is a $^*$-homomorphism $\psi:C_0(Y) \to M_m$ satisfying \eqref{eq:TracesClose}.

Without loss of generality, let us assume that $\mathcal H$ consists of positive contractions.
We may approximate the measure corresponding to $\tau$ by a convex combination (with rational coefficients) of point-mass measures.
Consequently, we may find some $k$ and a $^*$-homomorphism $\psi_0:C_0(Y) \to M_k$ such that
\begin{equation} |\tau(f)-\tau_{M_k}(\psi_0(f))|<\gamma/2,\quad f \in \mathcal H. \end{equation}

Choose a natural number $m_0 \geq \frac{2k}\gamma$.

Given $m\geq m_0$, write $m=kq+r$ where $0\leq r<k$.
Then define $\psi:=(\psi_0\otimes 1_q) \oplus \rho$, where $\rho:C_0(Y) \to M_r$ is any $^*$-homomorphism.
For $f \in \mathcal H$, we have
\begin{equation}
\tau_{M_m}(\psi(f)) = \frac{kq}m \tau_{M_k}(\psi_0(f)) + \frac rm \tau_{M_r}(\rho(f)). 
\end{equation}
Since both $\tau_{M_r}(\rho(f))$ and $\tau_{M_k}(\psi_0(f))$ are between $0$ and $1$, their difference is at most $1$, and so
\begin{equation} |\tau_{M_m}(\psi(f))-\tau_{M_k}(\psi_0(f))|<\frac rm < \frac km <\frac \gamma2, \end{equation}
and thus 
\begin{align}
|\tau(f)-\tau_{M_m}(\psi(f))| &\leq 
|\tau(f)-\tau_{M_k}(\psi_0(f))| + |\tau_{M_k}(\psi_0(f))-\tau_{M_m}(\psi(f))| \\
&< \gamma, 
\notag
\end{align}
as required.
\end{proof}

Making use of Lin's uniqueness theorem (Theorem \ref{thm:LinReformulated}) and the above lemma, we get the following lifting theorem for $^*$-homomorphisms from contractible spaces.

\begin{theorem}\label{thm:getting homoms}
Let $(m_n)_{n\in\mathbb N}$ be a sequence of natural numbers tending to $\infty$ and set $B_\infty:=\prod_n M_{m_n}/\bigoplus_n M_{m_n}$.
Let $Y$ be a second countable locally compact Hausdorff space for which $C_0(Y)$ is contractible.\footnote{Equivalently, the one-point compactification of $Y$ is contractible.}
If $\varphi:C_0(Y) \to B_\infty$ is a $^*$-homomorphism such that $\tau\circ\varphi$ is a faithful trace, for all $\tau \in T_\infty(B_\infty)$, then $\varphi$ lifts to a $^*$-homomorphism $\psi:C_0(Y) \to \prod_n M_{m_n}$; that is,
\begin{equation} q\circ \psi = \varphi, \end{equation}
where $q:\prod_n M_{m_n} \to B_\infty$ is the quotient map.
\end{theorem}

\begin{proof}
We shall construct a $^*$-homomorphism $C_0(Y) \to \prod_n M_{m_n}$ which agrees with $\varphi$ on all traces and then appeal to Lin's Uniqueness Theorem to say that $\varphi$ is unitarily equivalent to this map (composed with $q$).

By the Choi--Effros Lifting Theorem, lift $\varphi$ to a sequence of cpc maps $\varphi_n:C_0(Y) \to M_{m_n}$.

Let $(f_k)_{k\in\mathbb N}$ be a dense sequence in $C_0(Y)$.
For each $k$, let $m_0^{(k)}$ be the $m_0$ given by Lemma \ref{lem:TracesClose} for $\mathcal H:=\{f_1,\dots,f_k\}$ and $\gamma:=\frac1k$.
Let $n_k$ be such that $m_n\geq m_0^{(k)}$ for all $n\geq n_k$.
We may, in addition, arrange that $n_1<n_2<\cdots$.

By Lemma \ref{lem:TracesClose}, we may find a $^*$-homomorphism $\hat\psi_n:C_0(Y) \to M_n$, for each $n_k \leq n < n_{k+1}$, such that
\begin{equation} |\tau_{M_{m_n}}(\varphi_n(f_i))-\tau_{M_{m_n}}(\hat\psi_n(f_i))|<\frac1k,\quad i=1,\dots,k. \end{equation}
Defining $\hat\psi_1,\dots,\hat\psi_{n_1-1}$ arbitrarily, this yields a $^*$-homomorphism $\hat\psi:=(\hat\psi_n)_{n\in\mathbb N}:C_0(Y) \to \prod_n M_{m_n}$.
Since
\begin{equation} \lim_{n\to\infty} |\tau_{M_{m_n}}(\varphi_n(f_i))-\tau_{M_{m_n}}(\hat\psi_n(f_i))| = 0 \end{equation}
for all $i$, it follows that $\tau \circ \varphi = \tau\circ q\circ \hat\psi$ for all $\tau \in T_\infty(B_\infty)$.

Since $C_0(Y)$ is contractible, $\underline{K}(C_0(Y))=0$. 
Consequently, Theorem \ref{thm:LinReformulated} applies and tells us that there is a unitary $u \in B_\infty$ such that $\varphi=Ad_u \circ (q \circ \hat\psi)$.

Lifting $u$ to a unitary $v \in \prod_n M_{m_n}$\footnote{It is always possible to lift unitaries from sequence algebras, due to stability of the relation of being a unitary (see \cite[Corollary 2.22]{Blackadar}).} and setting $\psi:=Ad_v \circ \hat\psi$, we are done.
\end{proof}

The following corollary encapsulates how we wish to use the above theorem, by turning it into an order zero extension theorem, via the correspondence between order zero maps and $^*$-homomorphisms from a cone (see \cite[Corollary 3.1]{WZCPC0}).

In the following we make use of order zero functional calculus (see \cite[Corollary 3.2]{WZCPC0}).

\begin{corollary}\label{cor:Order0Extension}
Let $(m_n)_{n\in\mathbb N}$ be a sequence of natural numbers tending to $\infty$ and set $B_\infty:=\prod_n M_{m_n}/\bigoplus_n M_{m_n}$.
Let $X$ be a compact metrizable space, and let $\varphi:C(X) \to B_\infty$ be a cpc order zero map satisfying the following faithfulness condition:
for every trace $\tau \in T_\infty(B_\infty)$, every nonzero $f \in C_0((0,1])_+$, and every nonzero $g \in C(X)$,
\begin{equation} \tau(f(\varphi)(g))>0. \end{equation}
Then $\varphi$ extends to a cpc order zero map $\tilde\varphi:B(X) \to B_\infty$, where $B(X) \supset C(X)$ is the C*-algebra of bounded Borel functions on $X$.
\end{corollary}

\begin{proof}
By \cite[Corollary 3.1]{WZCPC0}, there is a $^*$-homomorphism $\hat\varphi:C_0((0,1]) \otimes C(X) \to B_\infty$ given by $\hat\varphi(f\otimes g) = f(\varphi)(g)$.
Since every nonzero positive element of $C_0((0,1])\otimes C(X)$ dominates one of the form $f\otimes g$ (with $f,g\geq 0$ both nonzero), we see that the faithfulness hypothesis implies that $\tau\circ\hat\varphi$ is a faithful trace on $C_0((0,1])\otimes C(X)$, for every $\tau \in T_\infty(B_\infty)$.
Therefore by Theorem \ref{thm:getting homoms}, $\hat\varphi$ lifts to a $^*$-homomorphism $\hat\psi=(\hat\psi_n)_{n\in\mathbb N}:C_0((0,1])\otimes C(X) \to \prod_n M_{m_n}$.
Since each $M_{m_n}$ is a von Neumann algebra, each $\hat\psi_n$ extends to a $^*$-homomorphism $\tilde\psi_n: C_0((0,1])\otimes B(X) \to M_{m_n}$ (in fact, we can even extend $\hat\psi_n$ to the larger algebra $B((0,1]\times X)$), and putting these $\tilde\psi_n$ together yields an extension of $\hat\psi$ to $\tilde\psi:C_0((0,1])\otimes B(X) \to \prod_n M_{m_n}$.
Using the converse direction in \cite[Corollary 3.1]{WZCPC0}, the $^*$-homomorphism $q\circ \tilde\psi:C_0((0,1])\otimes B(X) \to B_\infty$ (where $q:\prod_n M_{m_n} \to B_\infty$ is the quotient map) gives rise to a cpc order zero map $\tilde\varphi:B(X) \to B_\infty$.
Since $\tilde\psi$ extends $\hat\psi$, we see that $\tilde\varphi$ extends $\varphi$.
\end{proof}

\bigskip
\section{Certain quasicentral approximate units}
\label{sec:Quasicentral}
Recall, for a $C^*$-algebra $A$ and an ideal $I$, an approximate unit $(h_n)_{n\in\mathbb N}$ for $I$ is \emph{idempotent} if $h_{n+1}h_n=h_n$ for all $n \in \mathbb{N}$ and is \emph{quasicentral} (for ($I,A$)) if $\|h_n a-ah_n\| \to 0$ for all $a \in A$.
Assuming $A$ is separable, a quasicentral approximate unit always exists (\cite[Theorem 1]{Arvext}).

Given an idempotent approximate unit $(h_n)_{n\in\mathbb N}$ for $\mathcal K$, it follows that each $h_n$ must have finite rank, and therefore, the hereditary subalgebra $(h_{n+1}- h_n)B(\mathcal{H})(h_{n+1}-h_n)$ can be identified with a matrix algebra, $M_{m_n}$.
We make this identification in the following.

\begin{lemma}\label{lem:measures}
Let $X$ be a compact metrizable space and let $\tau$ be a trace on $C_0((0,1] \times X)$.
Then there exists a trivial extension $\T_\psi$ associated to a $^*$-homomorphism $\psi: C(X) \to B(\mathcal{H})$ (that is injective modulo $\K$) and an idempotent quasicentral approximate unit $(h_n)_\mathbb{N}$ for $(\mathcal K,\T_{\psi})$, such that the cpc maps
 \begin{equation}\psi_n: C((0,1] \times X) \to (h_{n+1}- h_n)B(\mathcal{H})(h_{n+1}-h_n)\cong M_{m_n}\end{equation}
 defined by 
\begin{equation}
\label{eq:measures}
\psi_n(\mathrm{id}_{(0,1]}^k \otimes g) := (h_{n+1}-h_n)^{\frac{k}{2}}(\psi)(g)(h_{n+1}-h_n)^{\frac{k}{2}}\end{equation}
(for $k \in \mathbb N, g \in C(X)$)
 satisfy $\lim_{n \to \infty} \tau_{M_{m_n}} \psi_n(f) = \tau(f)$ for all $f \in C_0((0,1]\times X)$. 
\end{lemma}

\begin{proof}
Let $(f_n)_{n\in\mathbb N}$ be a dense sequence in the set of positive contractions in $C_0((0,1]\times X)$.
For each $n$, by Lemma \ref{lem:TracesClose} (with $\mathcal H:=\{f_1,\dots,f_n\}$), choose some $r_n$ and some $^*$-homomorphism $\rho_n:C_0((0,1]\times X) \to M_{r_n}$ such that
\begin{equation} |\tau(f_i)-\tau_{M_{r_n}}(\rho_n(f_i))|<\frac1n,\quad i=1,\dots,n. \end{equation}
We may assume that $\frac{r_n}{r_{n+1}} \to 0$.

Without loss of generality (by decomposing into irreducible representations), $\rho_n$ has the form
\begin{equation} \rho_n(f)=\mathrm{diag}(f(t^{(n)}_1,x^{(n)}_1),\dots,f(t^{(n)}_{r_n},x^{(n)}_{r_n})) \end{equation}
for some $t^{(n)}_1,\dots,t^{(n)}_{r_n} \in (0,1]$ and some $x^{(n)}_1,\dots,x^{(n)}_{r_n} \in X$.
Without loss of generality (by possibly adding more points sporadically), we may also assume that for every $n_0$, $\{x^{(n)}_i: n \geq n_0\}$ is dense in $X$.

Define a $^*$-homomorphism $\psi:C(X) \to B(H)$ by 
\begin{equation} \psi(g):=\mathrm{diag}(g(x^{(1)}_1),\dots,g(x^{(1)}_{r_1}),g(x^{(2)}_1),\dots,g(x^{(2)}_{r_2}),\dots), \end{equation}
and for each $n$, define
\begin{equation}
h_n := \mathrm{diag}(\overbrace{1,\ \ \dots\ \ ,1}^{r_1+\cdots+r_{n-1}\text{ times}},t^{(n)}_1,\dots,t^{(n)}_{r_n},0,\dots) \in \mathcal K.
\end{equation}
Since $\{x^{(n)}_i: n \geq n_0\}$ is dense for every $n_0$, we can see that if $g\neq 0$ then $\psi(g)\not\in\K$; that is, $\psi$ is injective modulo $\mathcal K$, and thus $\psi$ defines a trivial extension $\mathcal T_\psi$.

Evidently, $(h_n)_{n\in\mathbb N}$ is an idempotent approximate identity for $\mathcal K$, and it commutes with the image of $\psi$ and is therefore quasicentral for $(\mathcal K, T_\psi)$.

Since each $h_n$ commutes with the image of $\psi$, we see that we have $\psi_n(f\otimes g)=f(h_{n+1}-h_n)\psi(g)$ for all $f\in C_0((0,1])$ and $g \in C(X)$.

We have
\begin{equation} h_{n+1}-h_n= \mathrm{diag}(\overbrace{0,\ \ \dots\ \ ,0}^{r_1+\cdots+r_{n-1}\text{ times}},1-t^{(n)}_1,\dots,1-t^{(n)}_{r_n},t^{(n+1)}_1,\dots,t^{(n+1)}_{r_{n+1}},0,\dots), \end{equation}
so its hereditary subalgebra is a copy of $M_{m_n}$ where $m_n=r_n+r_{n+1}$.
Making this identification, and defining $\psi_n$ as in \eqref{eq:measures}, we have
\begin{align}
\psi_n(f\otimes g)
&= \mathrm{diag}(f(1-t^{(n)}_1)g(x^{(n)}_1),\dots, f(1-t^{(n)}_{r_n})g(x^{(n)}_{r_n}), \\
&\qquad f(t^{(n+1)}_1)g(x^{(n+1)}_1),\dots, f(t^{(n+1)}_{r_{n+1}})g(x^{(n+1)}_{r_{n+1}})). 
\notag
\end{align}
Thus in particular, the bottom-right $r_{n+1}$ entries of $\psi_n$ exactly constitute the map $\rho_{n+1}$, and so for a positive contraction $f \in C_0((0,1]\times X)$,
\begin{equation} |\tau_{M_{m_n}}(\psi_n(f))-\tau_{M_{r_n}}(\rho_{n+1}(f))|\leq 2\frac{r_n}{r_n+r_{n+1}}. \end{equation}
Hence, for $i\leq n$,
\begin{align}
|\tau_{M_{m_n}}(\psi_n(f_i))-\tau(f_i)|
&\leq |\tau_{M_{m_n}}(\psi_n(f_i))-\tau_{M_{r_n}}(\rho_{n+1}(f_i))| \\
\notag
&\qquad + |\tau_{M_{r_n}}(\rho_{n+1}(f_i))-\tau(f_i))| \\
\notag
&\leq 2\frac{r_n}{r_n+r_{n+1}} + \frac1{n+1} \to 0.
\notag
\end{align}
It follows, by taking linear combinations and using density, that $\lim_{n \to \infty} \tau_{M_{m_n}} \psi_n(f) = \tau(f)$ for all $f \in C_0((0,1]\times X)$. 
\end{proof}

We now use the previous lemma to get a quasicentral approximate unit with tracial properties that we want, inside the direct sum of an arbitrary extension with a trivial extension.

\begin{corollary}\label{cor:extension-sums}
Let $X$ be a compact metrizable space, let $\T_\phi$ be an extension of $C(X)$ by $\K$, and let $\tau \in T(C_0((0,1]\times X))$.
Then there exists a trivial extension $\T_\psi$ and a quasicentral approximate unit $(h_n)_{n\in\mathbb N}$ for $(\mathcal K, \T_{\phi\oplus \psi})$ such that the cpc maps
\begin{equation}(\phi\oplus \psi)_n:C((0,1] \times X) \to (h_{n+1}-h_n)B(\Hil \oplus \Hil)(h_{n+1}-h_n)\cong M_{m_n}\end{equation}
defined by
 \begin{equation}(\phi\oplus \psi)_n(\mathrm{id}_{(0,1]}^k \otimes g) = (h_{n+1}-h_n)^{\frac{k}{2}}(\phi)(g)(h_{n+1}-h_n)^{\frac{k}{2}}\end{equation}
 satisfy $\lim_{n \to \infty} \tau_{M_{m_n}} (\phi\oplus \psi)_n(f) = \tau(f)$ for all $f \in C_0((0,1]\times X)$. 
\end{corollary}

\begin{proof}
Let $\T_\psi$ be a trivial extension and let $(h^{(\psi)}_n)_{n\in\mathbb N}$ be an idempotent quasicentral approximate unit for $(\K,\T_{\psi})$ as given by Lemma \ref{lem:measures}.
Let $(h^{(\phi)}_n)_{n\in\mathbb N}$ be any idempotent quasicentral approximate unit for $(\K,\T_{\phi})$ (by \cite[Lemma 7.3.1]{BO}, for example).
Let $\psi_n:C_0((0,1]\times X) \to M_{m_n^{(\psi)}}$ and $\phi_n:C_0((0,1]\times X) \to M_{m_n^{(\phi)}}$ be the maps corresponding to these extensions and quasicentral approximate units (as in \eqref{eq:measures}).
By possibly replacing $(h^{(\psi)}_n)_{n\in\mathbb N}$ by a subsequence, we may assume that
\begin{equation} \label{eq:RankQuotientAssumption}
 \frac{\mathrm{rank}(h^{(\phi)}_n)}{\mathrm{rank}(h^{(\psi)}_n)-\mathrm{rank}(h^{(\psi)}_{n-1})} \to 0. \end{equation}
Define
\begin{equation} h_n:=h_n^{(\phi)} \oplus h_n^{(\psi)}, \end{equation}
which forms an idempotent quasicentral approximate unit for $\T_{\phi\oplus \psi}$.
Moreover, we see that $m_n=m_n^{(\phi)}+m_n^{(\psi)}$ and $(\phi\oplus\psi)_n=\phi_n\oplus\psi_n$.
It follows from \eqref{eq:RankQuotientAssumption} that $\frac{m_n^{(\phi)}}{m_n^{(\psi)}}\to 0$, and so
\begin{equation} \lim_{n\to\infty} \tau_{M_{m_n}}((\phi\oplus\psi)_n(f))=\lim_{n\to\infty} \tau_{M_{m_n^{(\psi)}}}(\psi_n(f))=\tau(f), \end{equation}
for all $f \in C_0((0,1]\times X)$.
\end{proof}

\bigskip
\section{Proof of the theorem}
\label{sec:Proof}

A piece of the completely positive approximations we will form come from partitions of unity for the quotient $C(X)$, based on a ``coloured'' open cover, as in the following.
The ``colouring'' idea is that we can assign to each set a colour from a palette of size $d+1$ (which we do formally using subscripts from $0$ to $d$), so that any two sets of the same colour are disjoint.
In the following, we also introduce a Borel cover compatible with the open sets corresponding to the first colour.

\begin{lemma}\label{lem:POU}
Let $X$ be a compact metric space of dimension $d$ and let $\delta>0$.
Then there exist a finite open cover $\{U^{(i)}_j: i=0,\dots,d; j=1,\dots,r(i)\}$ and a finite Borel cover $\{Y_j:j=1,\dots,r(0)\}$ of $X$ such that:
\begin{enumerate}
\item Each set $U^{(i)}_j$ and $Y_j$ has diameter at most $\delta$;
\item For each $i$, the family $\{U^{(i)}_j: j=1,\dots,r(i)\}$ is pairwise disjoint, as is the family $\{Y_j: j=1,\dots,r(0)\}$;
\item For each $j$, $U^{(0)}_j \subseteq Y_j$.
\end{enumerate}
\end{lemma}

\begin{proof}
The existence of the first partition of unity, $\{U^{(i)}_j:i=0,\dots,d; j=1,\dots,s(i)\}$ is given by \cite[Proposition 1.5, (i)$\Rightarrow$(ii)]{KW}.
Set $r(0):=s(0)+s(1)+\cdots+s(d)$ and set $U^{(0)}_j=\emptyset$ for $j>s(0)$.
Also set $r(i):=s(i)$ for $i>0$.
Define $V_1,\dots,V_{r(0)}$ to be the sets $U^{(0)}_1,\dots,U^{(0)}_{s(0)},\dots,U^{(d)}_1,\dots,U^{(d)}_{s(d)}$ (in that order).
Then define 
\begin{equation} Y_j:=V_j\setminus \bigcup_{j'<j} V_{j'}, \end{equation}
for $j=1,\dots,r(0)$.
These are Borel, form a disjoint partition of unity, and satisfy $U^{(0)}_j\subseteq Y_j$ for all $j$ (in fact, for each $j$ we either have equality or $U^{(0)}_j=\emptyset$).
\end{proof}

We now have all the pieces required to prove our main theorem, which we restate for the convenience of the reader.

{\renewcommand*{\thelemma}{A}
\begin{theorem}
Let $X$ be a compact metrizable space and let $E$ be an extension of $C(X)$ by $\K$:
\begin{equation} 0 \to \K \to E \to C(X) \to 0. \end{equation}
Then $\nd E = \dim(X)$.
\end{theorem}
}

\begin{proof}[Proof of Theorem \ref{thm:Main}]
Let $d:=\dim(X)$; we will show $\nd E \leq d$ (the converse inequality follows from \cite[Proposition 2.3(iv)]{WZ}, since the quotient $C(X)$ has nuclear dimension $d$).
If $d=0$ then $X$ is totally disconnected, and so $E$ is a trivial extension by \cite[Theorem 1.15]{BDF77}.
Consequently, there is a quasicentral approximate unit for $(\K,E)$ consisting of projections, and so by \cite[Proposition 6.2]{KW}, $\nd E \leq \dr \T_{\phi} = \dim(X)=0$.
(Alternatively, one could make a direct argument that $E$ is AF.)
Hence in what follows we may assume $d>0$.

Let $\phi:C(X)\to B(\Hil)$ be a ucp map which is an injective $^*$-homomorphism modulo $\K$, such that $E\cong \T_\phi$.
Taking $\tau_0$ to be any faithful trace on $C_0((0,1]\times X)$, let $\T_\psi$ be the trivial extension and $h=(h_n)_{n\in \mathbb{N}}$ the idempotent quasicentral approximate unit for $\K$ in $\T_{\psi+\phi}$ obtained from Corollary \ref{cor:extension-sums}.
(When we apply Corollary \ref{cor:Order0Extension} later in the proof, we will make use of the fact that $\tau_0$ is faithful.)
Set $\Phi:= \psi\oplus \phi$ and note that $\mathcal{T}_{\Phi} \cong \mathcal{T}_{\phi}$ since $\T_\psi$ is a trivial extension (Theorem \ref{thm:BDF}).
The rest of the proof will show that $\nd(\T_{\Phi})\leq d$.

Define $\tilde{h}:=(h_{n+1})_{n\in \mathbb N}$ to be the same approximate unit as $h$ but with the index shifted by 1.
We may view both $h$ and $\tilde{h}$ as elements of $(\mathcal T_{\Phi})_\infty$ (as defined below).
Define 
\begin{align}
A_n&:=\overline{h_n\T_{\Phi}h_n}, \quad & A_\infty&:=\prod_n A_n/\bigoplus_n A_n \\
\notag
B_n&:= \overline{(h_{n+1}-h_n)\T_{\Phi}(h_{n+1}-h_n)}, \quad & B_\infty&:=\prod_n B_n/\bigoplus_n B_n \\
\notag
C_n&:= \overline{(1-h_{n+1})\T_{\Phi}(1-h_{n+1})},\quad &B_\infty&:=\prod_n C_n/\bigoplus_n C_n, \\
\notag
&& (\T_\Phi)_\infty&:=\prod_n \T_\Phi/\bigoplus_n \T_\Phi. 
\end{align}
The algebras $A_\infty,B_\infty$, and $C_\infty$ give a sort of overlapping decomposition of $(T_\Phi)_\infty$.
Define corresponding cpc maps 
\begin{align}
\alpha:&\T_\Phi\to A_\infty,\quad & \alpha(x)&:=(h_n^{\frac{1}{2}}xh_n^{\frac{1}{2}})_{n\in \mathbb{N}} \\
\notag
\beta:&C(X)\to B_\infty,\quad & \beta(f)&:=((h_{n+1}-h_n)^{\frac{1}{2}}\Phi(f)(h_{n+1}-h_n)^{\frac{1}{2}})_{n\in \mathbb{N}} \\
\notag
\gamma:&C(X)\to C_\infty,\quad &\gamma(f)&:=((1-h_{n+1})^{\frac{1}{2}}\Phi(f)(1-h_{n+1})^{\frac{1}{2}})_{n\in \mathbb{N}}.
\end{align}
Define $\iota_\infty:\T_\Phi \to (\T_\Phi)_\infty$ to be the diagonal inclusion.
Let $\iota_A:A_\infty \to (\T_\Phi)_\infty$, $\iota_B:B_\infty\to (\T_\Phi)_\infty$, and $\iota_C:C_\infty \to (\T_\Phi)_\infty$ be the inclusions coming from $A_n\subseteq \T_\Phi$, etc.
Finally, define $\pi:\T_\Phi \to C(X)$ to be the quotient map.

Using that $(h_n)_{n\in\mathbb N}$ is a quasicentral approximate unit, we see that $h,\tilde h$ commute with $\iota_\infty(\T_\Phi)$, and that for $x \in \T_\Phi$,
\begin{align}
\label{eq:abcNiceForm}
\iota_A \circ \alpha(x) &= h\iota_\infty(x), \\
\notag
\iota_B \circ \beta \circ \pi(x) &= (\tilde h-h)\iota_\infty(x), \\
\notag
\iota_C \circ \gamma \circ \pi(x) &= (1-\tilde h)\iota_\infty(x).
\end{align}
In particular, since the right-hand sides are all order zero functions of $x$, and since $\iota_A,\iota_B,\iota_C$ are injective while $\pi$ is surjective, it follows that $\alpha,\beta,\gamma$ are each cpc order zero maps.
Also, it follows that
\begin{equation} 
\label{eq:abcDecomp}
 \iota_\infty(x)=\iota_A\circ\alpha(x) + \iota_B\circ\beta\circ\pi(x)+\iota_C\circ \gamma\circ\pi(x). \end{equation}

We will fill in the maps in the following diagram (a variant on diagram in the proof of the Toeplitz case by Brake and Winter \cite{BW}), and its various paths will together form our approximation:
\begin{equation}\begin{tikzcd}
& F_k \arrow[d, dashed, shift left, "\varphi_k"] \arrow[r, dashed, "\rho_k"] & B(X) \arrow[r, "\tilde{\beta}"] & B_\infty \arrow[dr, "\iota_B"] &  \\
\T_\Phi \arrow[r,"\pi"] \arrow[drrr, "\alpha"'] & C(X) \arrow[u, dashed, shift left, "\psi_k"] \arrow[ur, hookrightarrow] \arrow[urr, "\beta"'] \arrow[rr, "\gamma"'] &  & C_\infty \arrow[r, "\iota_C"'] & (\T_\Phi)_\infty \\
& & & A_\infty \arrow[ur, "\iota_A"'] & 
\end{tikzcd}\end{equation}

By the choice of $(h_n)_{n\in\mathbb N}$ in Corollary \ref{cor:extension-sums}, the cpc order zero map $\beta$ satisfies $\tau(f(\beta)(g))=\tau_0(f\circ g)$, for any $f \in C_0((0,1])$ and any $g \in C(X)$.
In particular, since $\tau_0$ is faithful on $C_0((0,1]\times X)$, $\beta$ satisfies the faithfulness condition in Corollary \ref{cor:Order0Extension}, and so it extends to an order zero function $\tilde\beta:B(X) \to B_\infty$ (recalling that $B(X)$ denotes the bounded Borel functions on $X$).

For each $k$, choose an open partition of unity $\{U^{(i)}_{k,j}:i=0,\dots,d,\ j=1,\dots,r(k,i)\}$ and a Borel partition of unity $\{Y_{k,j}:j=1,\dots,r(k,0)\}$ satisfying the conditions of Lemma \ref{lem:POU}, where the bound on the radius is a quantity $\delta \to 0$.
Using the open cover, we shall define maps $\psi_k,\varphi_k$ using the idea of the proof that $\nd C(X)\leq d$ (see \cite[Proposition 2.19]{WCD1}).

For each $i,j,k$, choose any $x^{(i)}_{k,j}$ in $U^{(i)}_{k,j}$ (or in $Y_{k,j}$ if $i=0$ and $U^{(i)}_{k,j}=\emptyset$).
Also, select a partition of unity $(g^{(i)}_{k,j})_{i,j}$ in $C(X)$ subordinate to to the open cover $(U^{(i)}_{k,j})_{i,j}$.
Define $\psi_k=(\psi_k^{(i)})_{i=0}^d:C(X) \to \mathbb C^{\oplus r(k,0)}\oplus \cdots \oplus \mathbb C^{\oplus r(k,d)}$ by
\begin{equation} \psi_k^{(i)}(f):= (f(x^{(i)}_{k,1}),\dots,f(x^{(i)}_{k,r(k,i)})) \end{equation}
(a $^*$-homomorphism), $\varphi_k^{(i)}:\mathbb C^{\oplus r(k,i)} \to C(X)$ by
\begin{equation} \varphi^{(i)}_k(\lambda^{(i)}_1,\dots,\lambda^{(i)}_{r(k,i)}):=\sum_{j=1}^{r(k,i)} \lambda^{(i)}_j g^{(i)}_{k,j} \end{equation}
(a cpc order zero map), and $\rho_k:\mathbb C^{\oplus r(k,0)} \to B(X)$ by
\begin{equation} \rho_k(\lambda^{(0)}_1,\dots,\lambda^{(0)}_{r(k,0)}):=\sum_{j=1}^{r(k,0)} \lambda^{(i)}_j \chi_{Y_{k,j}} \end{equation}
(a $^*$-homomorphism).
Set 
\begin{equation} \varphi_k:=\sum_{i=0}^d \varphi_k^{(i)}:\mathbb C^{\oplus r(k,0)}\oplus \cdots \oplus \mathbb C^{\oplus r(k,d)} \to C(X). \end{equation}
Then 
\begin{equation}
\label{eq:dimXApproximation}
 \varphi_k\circ\psi_k(f) \to f \quad \text{and} \quad \rho_k \circ \psi_k^{(0)}(f) \to f,\quad f \in C(X)
\end{equation} (converging in norm, see the proof of \cite[Proposition 2.19]{WCD1}).

In particular, for $f \in C(X)$,
\begin{equation}
\label{eq:rhobeta}
\tilde\beta\circ\rho_k\circ\psi_k^{(0)}(f) \to \tilde\beta(f) = \beta(f).
\end{equation}

Let us now show that
\begin{equation}
\label{eq:rhoplusphi0}
 (\iota_B\circ \tilde\beta\circ\rho_k + \iota_C\circ \gamma \circ \varphi_k^{(0)}):\mathbb C^{\oplus r(k,0)} \to (\T_\Phi)_\infty \end{equation}
is an order zero map, for all $k$.
Since both $\rho_k$ and $\varphi_k^{(0)}$ are cpc order zero (and the other components are $^*$-homomorphisms), we only need to show that for $j\neq j'$,
\begin{equation} \big(\iota_B\circ \tilde\beta\circ\rho_k(e_j)\big)\, \big(\iota_C\circ \gamma \circ \varphi_k^{(0)}(e_{j'})\big) = 0, \end{equation}
where $e_1,\dots,e_{r(k,0)}$ are the canonical basis elements of $\mathbb C^{\oplus r(k,0)}$.

Pick a sequence of functions $(d_m)_{m\in \mathbb{N}}$ supported on $U^{(0)}_{k,j'}$ such that $d_m g^{(0)}_{k,j'} \to g^{(0)}_{k,j'}$.
Since
\begin{equation} U^{(0)}_{k,j'} \subseteq Y_{k,j'} \subseteq X \setminus Y_{k,j}, \end{equation}
it follows that $\chi_{Y_{k,j}}\leq 1-d_m$, and thus
\begin{eqnarray}
&&\hspace*{-2em} \| \iota_B\circ \tilde\beta\circ\rho_k(e_j)\cdot \iota_C\circ \gamma \circ \varphi_k^{(0)}(e_{j'})\|  \\
\notag
&=& \|\iota_B(\tilde\beta(\chi_{Y_{k,j}}))\iota_C(\gamma(g^{(i)}_{k,j'}))\| \\
\notag
&\leq& \|\iota_B(\beta(1-d_m))\iota_C(\gamma(g^{(i)}_{k,j'}))\| \\
\notag
&\stackrel{\eqref{eq:abcNiceForm}}=& \|(\tilde h-h)\iota_\infty(\Phi(1-d_m))(1-\tilde h)\iota_\infty(\Phi(g^{(i)}_{k,j'}))\| \to 0.
\end{eqnarray}
The map \eqref{eq:rhoplusphi0} is cp since it is the sum of cp maps, and it is contractive since
\begin{eqnarray}
(\iota_B\circ \tilde\beta\circ\rho_k + \iota_C\circ \gamma \circ \varphi_k^{(0)})(1_{\mathbb C^{\oplus r(k,0)}}) 
\hspace*{-2mm} &\leq& \hspace*{-2mm} \iota_B(\beta(1_{C(X)})+\iota_C(\gamma(1_{C(X)})) \\
\notag
&\stackrel{\eqref{eq:abcNiceForm}}=& \hspace*{-2mm} (\tilde h - h) + (1-\tilde h) = 1-h\leq 1.
\end{eqnarray}

Having established that the function \eqref{eq:rhoplusphi0} is cpc order zero, we will now lift it; this will form one of the upward maps in the definition of nuclear dimension.
Using order zero lifting (\cite[Remark 2.4]{KW}), we may lift the function \eqref{eq:rhoplusphi0} to a cpc order zero map 
\begin{equation} \hat\varphi_k^{(0)}=(\hat\varphi_{k,n}^{(0)})_n:\mathbb C^{\oplus r(k,0)} \to \prod_{n\in\mathbb N} \T_\Phi. \end{equation}
Likewise, we may lift each $\gamma\circ\varphi_k^{(i)}:\mathbb C^{\oplus r(k,i)} \to \prod_n C_n$ to a cpc order zero map
\begin{equation} \hat\varphi_k^{(i)}=(\hat\varphi_{k,n}^{(i)})_n:\mathbb C^{\oplus r(k,i)} \to \prod_{n\in\mathbb N} C_n \end{equation}

Write $\alpha_n:\T_\Phi \to A_n$ for the map $\alpha_n(x):=h_n^{1/2} x h_n^{1/2}$ and $\iota_{A_n}:A_n \to \T_\Phi$ for the inclusion.
Note that $\iota_{A_n}+\hat\varphi_{k,n}^{(1)}:A_n \oplus \mathbb C^{\oplus r(k,1)} \to \T_\Phi$ is cpc order zero, since it is the sum of two cpc order zero maps with orthogonal ranges.

Altogether, we obtain an approximate factorization through $\mathbb C^{\oplus r(k,0)} \oplus (A_n \oplus \mathbb C^{\oplus r(k,1)}) \oplus \mathbb C^{\oplus r(k,1)} \oplus \cdots \oplus \mathbb C^{\oplus r(k,d)}$, given by the cpc maps
\begin{align}
\psi_k^{(0)}\circ \pi&:\T_\Phi \to \mathbb C^{\oplus r(k,0)}, \\
\notag
(\alpha_n\oplus \psi_k^{(1)}\circ \pi)&:\T_\Phi \to (A_n \oplus \mathbb C^{\oplus r(k,1)}), \\
\notag
\psi_k^{(i)}\circ \pi&:\T_\Phi \to \mathbb C^{\oplus r(k,i)}, \quad i=2,\dots,d
\end{align}
(where $\pi:\T_\Phi \to C(X)$ is the quotient map) together with the cpc order zero maps
\begin{align}
\hat\varphi_{k,n}^{(0)}&:\mathbb C^{\oplus r(k,0)} \to \T_\Phi, \\
\notag
(\iota_{A_n}+\hat\varphi_{k,n}^{(1)})&:A_n \oplus \mathbb C^{\oplus r(k,1)} \to \T_\Phi, \\
\notag
\hat\varphi_{k,n}^{(i)}&:\mathbb C^{\oplus r(k,i)} \to C_n \subseteq \T_\Phi, \quad i=2,\dots,d.
\end{align}
Let us write $\Theta_{k,n}:\T_\Phi \to \T_\Phi$ for the entire composition; that is,
\begin{equation} \Theta_{k,n}(x)=\alpha_n(x)+\sum_{i=0}^d \hat\varphi_{k,n}^{(i)}(\psi_k^{(i)}(\pi(x))). \end{equation}
We are done once we show that $\Theta_{k,n}$ converges in point-norm to $\mathrm{id}_{\T_\Phi}$.

Fixing $k$, by working in $(\T_\Phi)_\infty$, and using the maps that the $\hat\varphi_{k,n}^{(i)}$ lift, we have
\begin{align}
&\hspace*{-4mm}\limsup_{n\to\infty} \|\Theta_{k,n}(x)-x\| \\
\notag
&=\|\iota_A(\alpha(x))+\iota_B(\tilde\beta\circ\rho_k(\psi_k^{(0)}(\pi(x)))) + 
\sum_{i=0}^d \iota_C\circ\gamma\circ\varphi_k^{(i)}\circ\psi_k^{(i)}\circ\pi(x)\| \\
&=\|\iota_A\circ\alpha(x)+\iota_B\circ\tilde\beta\circ\rho_k\circ\psi_k^{(0)}\circ\pi(x) + \iota_C\circ\gamma\circ\varphi_k\circ\psi_k\circ\pi(x)\|
\notag
\end{align}
Thus,
\begin{eqnarray}
&&\hspace*{-22mm} \limsup_{k\to\infty} \limsup_{n\to\infty} \|\Theta_{k,n}(x)-x\| \\
\notag
&\leq& \limsup_{k\to\infty} \|\iota_A\circ\alpha(x)+\iota_B\circ\tilde\beta\circ\rho_k\circ\psi_k^{(0)}\circ\pi(x) \\
\notag
&& \quad + \iota_C\circ\gamma\circ\varphi_k\circ\psi_k\circ\pi(x)\| \\
\notag
&\stackrel{\eqref{eq:dimXApproximation}\eqref{eq:rhobeta}}=& \|\iota_A\circ\alpha(x)+\iota_B\circ\beta\circ\pi(x) + \iota_C\circ\gamma\circ\pi(x)\| \\
\notag
&\stackrel{\eqref{eq:abcDecomp}}=&0,
\end{eqnarray}
as required.
\end{proof}

\end{document}